\documentclass[authoryear,11pt]{elsarticle}
\usepackage[top=2.6cm,bottom=2.8cm,left=3cm,right=3cm]{geometry}
\usepackage{graphicx}%
\usepackage{multirow}%
\usepackage{amsmath,amssymb,amsfonts}%
\usepackage{amsthm}%
\usepackage{mathrsfs}%
\usepackage[title]{appendix}%
\usepackage{xcolor}%
\usepackage{textcomp}%
\usepackage{manyfoot}%
\usepackage{booktabs}%
\usepackage{algorithm}%
\usepackage{adjustbox}
\usepackage{algpseudocode}%
\usepackage{listings}%
\usepackage{hyperref}
\usepackage{amsmath,amssymb,amsthm,mathtools}
\usepackage{url}
\usepackage{svg}
\usepackage{float}
\usepackage{multirow}
\usepackage{diagbox}

\usepackage{caption}
\usepackage{subcaption}

\newtheorem{definition}{Definition}
\newtheorem{proposition}{Proposition}

\newtheorem{example}{Example}
\newtheorem{theorem}{Theorem}
\newtheorem{corollary}{Corollary}

\newcommand{\R}{\mathbb{R}}

\newcommand{\pt}{\text{ }\forall\text{ }}
\newcommand{\tq}{\text{ }:\text{ }}
\newcommand{\N}{\mathbb{N}}

\newcommand{\io}{[0,1]}

\newcommand{\cR}{\mathcal{R}}
 \newcommand{\rmd}{\mathcal{R}_\delta}
 \newcommand{\rae}{\mathcal{R}_\varepsilon^\Phi}
\newcommand{\rag}{\mathcal{R}_\gamma^{w,\Phi}}
\newcommand{\rmg}{\mathcal{R}_\gamma^{w}}
\newcommand{\rowae}{\mathcal{R}_\varepsilon^{\Psi_\omega}}

\date{January 2025}
\journal{""}
\begin{document}
\begin{frontmatter}
		\title{Mutual Consensus and its Application in Minimum Cost Consensus Models.}
		\author[1]{Diego Garc\'ia-Zamora$^*$}
		\author[2]{Bapi Dutta}
        \author[3]{Luis Martinez}
		\address[1]{Department of Mathematics, Universidad de Ja{\'e}n,		Ja{\'e}n, Spain {dgzamora@ujaen.es}}
        \address[2]{Department of Computer Science, Universidad de Ja{\'e}n,		Ja{\'e}n, Spain {bdutta@ujaen.es}}
        \address[3]{Department of Computer Science, Universidad de Ja{\'e}n,		Ja{\'e}n, Spain {martin@ujaen.es}}

		\begin{abstract}
			This paper introduces the concept of {mutual consensus} as a novel non-compensatory consensus measure that accounts for the maximum disparity among opinions to ensure robust consensus evaluation. Incorporating this concept, several new Minimum Cost Consensus (MCC) models are proposed, and their properties are analyzed. To show their applicability, these mutual consensus-based MCC models are then considered in the context of the {OWA-MCC} model, which employs Ordered Weighted Averaging (OWA) operators for preference aggregation. Concretely, we include a linearized formulation under symmetry conditions as well as examples of the non-convexity of the feasible region in the general case.  Finally, mutual consensus is utilized to obtain approximate solutions for the OWA-MCC model, demonstrating its practical effectiveness and advancing the theoretical and applied dimensions of consensus modeling in group decision-making.
		\end{abstract}

		
	\end{frontmatter}
\section{Introduction}

Group Decision Making (GDM) provides a structured framework for collective decision scenarios in which multiple decision-makers (DMs) evaluate and prioritize alternatives based on their individual preferences \citep{salo1995interactive,herrera2005managing}. In such settings, disagreements among participants are common, especially when opinions are subjective or based on imprecise information \citep{guetzkow1954analysis}. To address this, Consensus Reaching Processes (CRPs) have been developed to iteratively guide the group toward a collectively agreed solution by suggesting changes to DMs’ opinions \citep{dong2018consensus,han2021robust}. In this context, it must be acknowledged that aiming for a full agreement is unfeasible in practice, especially in large-scale GDM (LSGDM) problems that may involve hundreds or thousands of DMs. This has motivated the adoption of the soft consensus paradigm, which promotes a more flexible interpretation of consensus \citep{kacprzyc1986}. Rather than requiring unanimity, soft consensus measures the degree of agreement after each round of interaction and concludes the process once a predefined threshold has been reached.

Among the different consensus models in the literature, Minimum Cost Consensus (MCC) models have gained prominence as a mathematical optimization-based approach that minimizes the cost required to modify individual opinions under predefined consensus thresholds \citep{benarieh2007}. Despite their utility, classical MCC models often rely on traditional consensus measures that primarily evaluate either the similarity of individual opinions to a collective judgment or the average mutual resemblance among individual opinions \citep{fzz-mcc}. These measures, while effective, may not fully capture scenarios where ensuring agreement across all pairs of opinions is critical. To address this gap, we introduce the concept of {mutual consensus}, a novel consensus measure that focuses on minimizing the maximum disagreement between any pair of opinions. This measure redefines the measurement of consensus, emphasizing robustness by ensuring that no extreme disparities remain in the group’s collective decision. 

Building upon this new measure, we introduce MCC models that incorporate mutual consensus as a central constraint. The mathematical properties of these models are then analyzed, offering insights into their behavior and potential to enhance decision-making processes. In particular, we provide alternative versions of these models with enhanced computational performance. Furthermore, we demonstrate that the mutual consensus measure is stronger than other classical consensus measures. Concretely, given a consensus threshold for these classical consensus measures, if the mutual consensus measure falls below this threshold, then the original consensus measure also falls below the same threshold.

We also present an application of mutual consensus in the context of {OWA-MCC} models \citep{zhang2011}, which leverages Ordered Weighted Averaging (OWA) operators to aggregate preferences in the classical MCC. We first provide a linearized version of the OWA-MCC model under symmetry conditions, enabling efficient computation in practical applications. However, we also identify and analyze examples where linearization is not computationally feasible, which lead to NP-hard problems that cannot be efficiently solved for large numbers of DMs. To address this issue, this manuscript explores the relationship between mutual consensus and classical OWA-based MCC models. By establishing theoretical connections, we obtain an efficient algorithm based on mutual consensus that allows computing an approximation to the OWA-MCC problem. Computationally, such an algorithm highlights itself because of its theoretical properties and computational efficiency, which provides a feasible alternative to classical Binary Linear Programming (BLP)-based strategies used for solving OWA-MCC \citep{zhang2013}.

The remainder of this manuscript is as follows. Section \ref{sec:prelim} introduces the basic notions to understand the proposal. In Section \ref{sec:MC}, we formally introduce the concept of {mutual consensus} and outline its theoretical foundations, highlighting its distinction from classical consensus measures. We also develop new MCC models based on mutual consensus and examine their mathematical properties. In Section \ref{sec:OWAMCC}, mutual consensus is integrated into the classical {OWA-MCC} model, with a detailed analysis of its linearization under symmetry conditions and examples illustrating the limits of the convexity of its feasible region. Section \ref{sec:relationship} explores the relationship between mutual consensus and classical OWA-based MCC models, providing a comparative analysis of their strengths and applications. Section \ref{sec:algorithm} demonstrates the practical applicability of mutual consensus through an approximation approach to solving the OWA-MCC model, which is verified in the computational analysis developed in Section \ref{sec:comput}. Finally, we present some concluding remarks and directions for future research in Section \ref{sec:conclusion}.

\section{Preliminaries}\label{sec:prelim}
 In this section, we introduce the basic concepts necessary to frame our proposal, including a review of GDM and consensus and the main elements of MCC models and their extensions under different approaches.

A GDM problem involves a set of DMs $E = \{e_1, e_2, \ldots, e_n\}$ and a set of alternatives $X = \{x_1, x_2, \ldots, x_m\}$. DMs express their preferences regarding the alternatives through preference structures, which can be represented as numerical vectors, pairwise comparison matrices, or fuzzy preference relations. The goal is to reach a group decision that may be acceptable to all DMs \citep{kacprzyc2010}. 

To ensure that all DMs feel that their opinions are reflected in the final decision, CRPs were proposed as iterative discussion processes in which a moderator leads the group towards an agreement \citep{ElSurvey}. Since a total agreement is nearly unfeasible when dealing with large groups of DMs, the notion of soft consensus was proposed as an alternative definition that considers instead degrees of consensus \citep{kacprzyc1986}. In this case, the CRP is assumed to be successful when the group consensus degree is high enough. The degree of acceptance of a solution to a GDM problem is evaluated through consensus measures. These measures can be classified into two main types \citep{paramcons}: whereas the first class evaluates the similarity between individual opinions and the group's collective opinion, the second category considers the similarity among individual opinions without referencing a collective opinion. Mathematical formulations for these measures include aggregation operators and distance functions that quantify the degree of consensus achieved \citep{fzz-mcc}.

MCC models, which were originally proposed by \cite{benarieh2007}, formulate GDM problems as mathematical programming problems. Their goal is to minimize the cost associated with the changes needed in DMs' opinions to reach a predefined consensus level. 
\begin{definition}[MCC model \citep{benarieh2007,zhang2011}] Let us consider $n$ decision-makers who provide their opinions $o_1,...,o_n\in\io$ using a numerical scale. Assume that $c_1,...,c_n$ are the costs of modifying their opinions and $\Phi:\io^n\to\io$ is an averaging aggregation operator used to form a collective opinion \citep{averaginghumberto}. Then, the MCC model to obtain the consensual opinions of $n$ DMs, $\{x_1, ..., x_n\}$ can be formulated as the following optimization problem:

\begin{equation}
    \label{MCC}\tag{MCC}
    \begin{split}
    \min_{x\in\io^n} & \hskip 5pt \sum_{k=1}^nc_k|x_k-{o}_k|\\
    \text{s.t.} & \hskip 5 pt\begin{cases}
    {g}=\Phi(x_1,...,x_n)\\
    |x_k-g|\leq \varepsilon,\; k=1,2,...,n\\
    \end{cases}
    \end{split}
    \end{equation}
where $\epsilon$ is the maximum allowable distance between group and individual, which is interpreted as a soft consensus threshold.
\end{definition}

Over time, various extensions of MCC models have been proposed to address specific characteristics of GDM problems. For instance, \cite{rodriguez2021CMCC} extended the original model to account for additional consensus measures as linear constraints, whereas 
 \cite{fzz-mcc} generalized MCC models to admit any preference structure and many kinds of consensus constraints. These extensions enable MCC models to be applied across a wide variety of contexts,  especially in LSGDM problems.

 \cite{zhang2011} studied the case in which the aggregation operator is the OWA operator.
\begin{definition}[OWA operator \citep{yager1988}]
Let $\omega\in\io^n$ be a weighting vector such that $\sum_{i=1}^m{\omega_i}=1$. The OWA Operator $\Psi_\omega:\io^n\to\io$ associated to $w$ is defined by:
\begin{equation*}	\Psi_\omega({x})=\sum_{k=1}^n{\omega_kx_{\sigma(k)}}\pt {x}\in\io^n
\end{equation*}
where $\sigma$ is a permutation of the $n$-tuple $(1,2,...,n)$ such that $x_{\sigma(1)}\geq x_{\sigma(2)}\geq ...\geq x_{\sigma(n)}$.
\end{definition}
Subsequent studies have further explored Zhang et al.'s model under OWA operators, including a Binary Linear Programming (BLP) formulation \citep{zhang2013}. More recently, MCC models featuring quadratic cost functions and OWA operators have been investigated by \cite{zhang2022consensus}. Additionally, \cite{qu2023robust} examined MCC models in which the cost associated with the opinion modifications is uncertain, applying a robust optimization approach.

While several researchers have integrated OWA operators into GDM problems \citep{WU2022, xie2024ordered}, the OWA operator induces a nonlinear constraint in the optimization model that depends on the order of the variables. This usually implies an elevated computational cost to solve even simple instances that involve fewer than 10 DMs. Until now, no study has solved the general OWA-MCC problem for a large number of DMs in a limited time to be applicable in real-life scenarios.

\section{Mutual Consensus and its applications to Minimum Cost Consensus models}\label{sec:MC}
This section presents the main notion of this manuscript, the mutual consensus measure. Afterwards, we integrate this notion into MCC models and derive some properties of the resulting framework.

\begin{definition}
    We define the mutual consensus measure as the mapping $\kappa:\io^n\to\io$ given by \begin{equation*}
        \kappa(x)=\max_{1\leq i<j\leq n}\{\vert x_i-x_j\vert\}\pt x\in\io^n.
    \end{equation*}
\end{definition}

Since the mutual consensus measure focuses on the maximum disagreement, it captures critical disparities, providing a robust alternative to classical measures that often rely on averages and may miss outliers. In this sense, this consensus measure is non-compensatory, which ensures that high agreement among some opinions does not offset significant disagreements elsewhere. It is also symmetric, meaning that the computed disparity remains unchanged regardless of the order of comparisons, and is easy to combine with a scale factor, making it applicable across normalized preference scales. Note that here we interpret the consensus measure as a distance, instead of a similarity, as in other classical consensus measures. Although it would be straightforward to restate this definition in such terms, we prefer to keep this version for the sake of simplicity, following the discussion provided in \cite{fzz-mcc}. In any case, the mapping $1-\kappa$ is a consensus measure in the sense of \cite{Beliakov2014}.

The mutual consensus measure serves as the foundation for a new MCC model.
\begin{definition}[MCMC model] Let us consider a GDM problem involving $n$ decision-makers who provide their opinions $o_1,...,o_n\in\io$ using a numerical scale. Assume that $c_1,...,c_n$ are the costs of modifying their opinions and $\delta\in\io$ is a consensus threshold. Then, the Minimum Cost with Mutual Consensus (MCMC) model 
can be formulated as follows:

    \begin{equation}\label{MCMC}
    \tag{MCMC}
    \begin{split}
    \min_{x\in\io^n} & \hskip 5pt \sum_{k=1}^nc_k|x_k-{o}_k|\\
    \text{s.t.} & \hskip 5 pt\begin{cases}
    |x_i-x_j|\leq \delta,\; 1\leq i<j\leq n.\\
    \end{cases}
    \end{split}
    \end{equation}
\end{definition}

Now, let us introduce some notation for MCC models. First, given an initial opinion $o=(o_1,...,o_n)\in\io^n$, we will denote the cost function associated with the cost vector $c=(c_1,...,c_n)\in\io^n, ||c||_1=1$  as $\xi_c^o:\io^n\to\io$, which is defined by 
\begin{equation*}
    \xi_c^o(x)=\sum_{k=1}^nc_k|x_k-{o}_k|.
\end{equation*}
Note that considering the cost as a normalized vector allows removing the effect of the scale from the model without modifying the results. At this stage, the previous MCMC model can be written as 
\begin{equation}
    \tag{MCMC}
    \min_{x\in\rmd} \hskip 5pt \xi_c^o(x)
\end{equation}
where  $\delta\in\io$ and $\rmd=\{x\in\io^n\tq \kappa(x)\leq\delta\}$.
    
Note that we do not need to consider an aggregation operator in the definition of MCMC as in MCC, but we have $\frac{n(n-1)}{2}$ constraints in the feasible region. 
The following result allows the simplification of the original MCMC model by reducing the number of constraints from $\mathcal{O}(n^2)$  to  $\mathcal{O}(n)$.
\begin{theorem}
    Let us fix $c\in\io^n$ and  $o\in\io^n$. Consider a permutation $\sigma$ that orders decreasingly the opinion $n$-tuple $o$. If $x^*\in\io^n$ is an optimal solution of 
    \begin{gather*}  
    \min_{x\in\io^n}\xi_{c_\sigma}^{o_\sigma}(x)\\
    \text{s.t.}\begin{cases}
    x_1\geq x_2\geq...\geq x_n\\
    x_1-x_n\leq\delta
    \end{cases}
    \end{gather*}
    then, $x^*_{\sigma^{-1}}$ is an optimal solution for the MCMC problem $\min_{x\in\rmd} \hskip 5pt \xi_c^o(x)$.
\end{theorem}
\begin{proof}
     Since $x^*$ is a solution of the first optimization problem, then $\kappa(x^*)=x_1-x_n\leq\delta$, and, after reordering the terms in $x^*$ according to the permutation $\sigma^{-1}$, they will still satisfy the same restriction. Consequently, $\kappa(x^*_{\sigma^{-1}})\leq\delta$ and $x^*_{\sigma^{-1}}\in\rmd$. Therefore, $x^*_{\sigma^{-1}}$ is a feasible solution to the second optimization problem. 
     
     Further, for any $y\in\io^n$:
    \begin{gather*}       \xi_c^o(y_{\sigma^{-1}})=\sum_{i=1}^nc_i|y_{\sigma^{-1}(i)}-o_i|=\sum_{i=1}^nc_{\sigma(i)}|y_{i}-o_{\sigma(i)}|=\xi_{c_\sigma}^{o_\sigma}(y).
    \end{gather*}
    Therefore, for $y\in\arg\min_{x\in\rmd} \hskip 5pt \xi_c^o(x)$, since $o_\sigma$ is decreasingly ordered and $y$ is a minimum cost solution to $\arg\min_{x\in\rmd}\xi_c^o(x)$, then $y_\sigma$ is also decreasingly ordered and thus is a feasible solution for the first optimization problem, from which $\xi_c^o(y)=\xi_{c_\sigma}^{o_\sigma}(y_\sigma)= \xi_{c_\sigma}^{o_\sigma}(x_{\sigma^{-1}}^*)$. Consequently, $x_{\sigma^{-1}}^*$ is also minimal for $\min_{x\in\rmd} \hskip 5pt \xi_c^o(x)$.
    
\end{proof}

At this stage, if we fix an averaging aggregation operator $\Phi:\io^n\to\io$, we can define the consensus measure $\kappa^\Phi:\io^n\to\io$ as $\kappa^\Phi(x)=\max_{i=1,...,n}\{|x_i-\Phi(x)|\}$ and the associated feasible region $\rae=\{x\in\io^n\tq \kappa^ \Phi(x)\leq\varepsilon\}$ for a certain $\varepsilon\in\io$. In such a case, the classical MCC model can be written as 
\begin{equation}
    \tag{MCC}
    \min_{x\in\rae} \hskip 5pt \xi_c^o(x).
\end{equation}
Whenever we consider the region  $\rae$, we implicitly assume an averaging aggregation operator $\Phi$, although we do not explicitly mention it. Let us show the first important property of mutual consensus, i.e., if everyone in the group agrees with each other (within a tolerance level $\alpha$), then each individual also agrees with the group opinion, formed by any averaging aggregation operator $\Phi$ (within a tolerance level $\alpha$).
\begin{proposition}\label{P:MC-Inside}
    Let us consider $\alpha\in\io$. Then $\mathcal{R}_\alpha\subseteq\mathcal{R}^\Phi_\alpha$. In other words, if $\kappa(x)\leq\alpha$, then $\kappa^\Phi(x)\leq\alpha$.
\end{proposition}
\begin{proof}
    Let us consider $x\in\mathcal{R}_\alpha$. since $\Phi$ is an averaging aggregation function, $\min{(x)}\leq \Phi(x)\leq \max{(x)}$. Therefore, for each $k=1,...,n$, the following inequality holds:
    \begin{equation*}
        \begin{split}
            |x_k-\Phi(x)|\leq\max\{|x_k-\max{(x)}|,|x_k-\min{(x)}|\}\leq \alpha.
        \end{split}
    \end{equation*}
   This imples that  $x\in\mathcal{R}_\alpha^\Phi$.
\end{proof}

In addition to the maximum distance to the average $\kappa^\Phi$, \cite{rodriguez2021CMCC} considered MCC models that included, as constraints, the consensus measures $\kappa^{w,\Phi},\kappa^{w}:\io^n\to\io$ defined as follows:
\begin{equation*}
    \begin{split}
        \kappa^{w,\Phi}(x)&=\sum_{k=1}^nw_k|x_k-\Phi(x)|\pt x\in\io^n,\\
        \kappa^{w}(x)&=\sum_{k=1}^{n-1}\sum_{l=k+1}^n\frac{w_k+w_l}{n-1}|x_k-x_l|\pt x\in\io^n, 
    \end{split}
\end{equation*}
where $w\in\io^n$ is a weighting vector, i.e. it satisfies $\sum_{k=1}^nw_k=1$.
Thus, given a parameter $\gamma\in\io$, we can define the associated feasible regions as
\begin{equation*}
    \begin{split}
           \rag&=\{x\in\io^n\tq   \kappa^{w,\Phi}(x)\leq \gamma\},\\
        \rmg&=\{x\in\io^n\tq \kappa^w(x)\leq \gamma\}.
    \end{split}.
\end{equation*}
These regions correspond to the Comprehensive MCC (CMCC) models presented by \cite{labella2020}, which, according to our notation, are as follows:
\begin{equation}\label{CMCC}
\tag{CMCC} 
    \begin{split}
        \min_{x\in\rae\cap\rag} \hskip 5pt \xi_c^o(x)\\
        \min_{x\in\rae\cap\rmg} \hskip 5pt \xi_c^o(x)
    \end{split}
\end{equation}
Let us show that the notion of mutual consensus indeed implies all these other consensus measures.
\begin{theorem}\label{T:MCMC-General}
    Let us consider $\alpha\in\io$. Then, 
    \begin{equation*}
        \mathcal{R}_\alpha\subseteq \mathcal{R}_\alpha^{\Phi}\cap\mathcal{R}_\alpha^{w,\Phi}\cap\mathcal{R}_\alpha^{w}.
    \end{equation*}
\end{theorem}
\begin{proof}
    From Proposition \ref{P:MC-Inside}, we know that  $ \mathcal{R}_\alpha\subseteq \mathcal{R}_\alpha^{\Phi}$. Now, let us pick $x\in \mathcal{R}_\alpha$. Since $\kappa(x)\leq\alpha$, and $\kappa^{\Phi}(x)\leq\alpha$,  we obtain:
    \begin{equation*}
        \begin{split}
            \kappa^{w,\Phi}(x)&=\sum_{k=1}^nw_k|x_k-\Phi(x)|\leq \sum_{k=1}^nw_k\alpha=\alpha\\
            \kappa^{w}(x)&=\sum_{k=1}^{n-1}\sum_{l=k+1}^n\frac{w_k+w_l}{n-1}|x_k-x_l|\leq\sum_{k=1}^{n-1}\sum_{l=k+1}^n\frac{w_k+w_l}{n-1}\alpha\\
            &=\frac{\alpha}{n-1}\sum_{k=1}^{n-1}\sum_{l=k+1}^n(w_k+w_l)=\frac{\alpha}{n-1}(\sum_{k=1}^{n-1}w_k(n-k)+\sum_{k=1}^{n-1}\sum_{l=k+1}^nw_l)\\
            &=\frac{\alpha}{n-1}(\sum_{k=1}^{n-1}w_k(n-k)+\sum_{k=2}^{n}(k-1)w_k)\\
            &=\frac{\alpha}{n-1}(w_1(n-1)+\sum_{k=2}^{n-1}w_k(n-k)+\sum_{k=2}^{n-1}(k-1)w_k+w_n(n-1))\\
            &=\frac{\alpha}{n-1}(\sum_{k=1}^{n}(n-1)w_k)=\alpha.
        \end{split}
    \end{equation*}
\end{proof}
The previous result enhances the most important property of mutual consensus: it is stronger than any other consensus measure in the literature in the sense that if a set of opinions satisfies the mutual consensus constraint, it will also satisfy the constraints determined by any other classical consensus measure. Furthermore, an optimal solution of MCMC presents other interesting properties.
\begin{proposition}\label{prop:MCMC_xvalues}
    Let $\delta\in[0,1]$, $c\in\io^n$ and $o\in\io^n$ and consider $x^*\in\operatorname{arg}\min_{x\in\mathcal{R}_\delta}\xi_c^o(x)$. Then, $x^*_i\in\{x^*_{\sigma(1)},x^*_{\sigma(n)},o_i\}$ for $i=1,...,n$, where $\sigma$ is the permutation that decreasingly orders $x^*$. 
\end{proposition}
\begin{proof}
    Let us fix $i=1,...,n$. We have three cases:
    \begin{itemize}
        \item $o_i\geq x^*_{\sigma(1)}$. Suppose that $x^*_i<x^*_{\sigma(1)}$. Then, we have that $|x^*_i-o_i|> |x^*_{\sigma(1)}-o_i|$, which means that the $n$-tuple $y\in\io^n$ given by 
        \begin{equation*}
            y_j=\begin{cases}
                x^*_j&\text{ if }j\neq i\\
                x^*_{\sigma(1)}&\text{ if }j=i
            \end{cases}
        \end{equation*}
        satisfies $y\in\mathcal{R}_\delta$ and $\xi_c^o(y)<\xi_c^o(x^*)$, which is a contradiction. Consequently, it must be $x^*_i=x^*_{\sigma(1)}$.
        \item $o_i\leq x_{\sigma(n)}$. Reasoning analogously to the previous case leads to $x^*_i=x^*_{\sigma(n)}$. 
        \item $o_i\in]x^*_{\sigma(n)},x^*_{\sigma(1)}[$. In this case, the choice $x^*_i=o_i$ is the one that minimizes the cost function.
    \end{itemize}
\end{proof}

In other words, if the opinion $o_i$ is changed, it must change to be either the minimum or the maximum of the new opinion vector $x^*$. Below, we include two other results that provide some insights regarding the connection of optimal solutions with the MCMC models when they only differ on the mutual consensus parameter.

\begin{proposition}\label{P:d1d2MCM}
    Let $\delta_1,\delta_2\in[0,1]$ such that $\delta_1<\delta_2$, $c\in\io^n$ and $o\in\io^n$ verifying $\max\{o_1,...,o_n\}-\min\{o_1,...,o_n\}>\delta_2$. Consider $x^1\in\operatorname{arg}\min_{y\in\mathcal{R}_{\delta_1}}\xi_c^o(y)$ and $x^2\in\operatorname{arg}\min_{y\in\mathcal{R}_{\delta_2}}\xi_c^o(y)$ satisfying 
    \begin{gather*}
        x^2_{\sigma^2(n)}\leq x^1_{\sigma^1(n)}\leq  x^1_{\sigma^1(1)}\leq x^2_{\sigma^2(1)}.
    \end{gather*} Then, for $j=1,...,n$ 
    \begin{gather*}
        x^1_{\sigma^1(1)}-x^1_j\leq x^2_{\sigma^2(1)}-x^2_j\\
         x^1_j-x^1_{\sigma^1(n)}\leq x^2_j -x^2_{\sigma^2(n)}.
    \end{gather*}
    where $\sigma^1$ is the permutation that decreasingly orders $x^1$ and $\sigma^2$ is the permutation that decreasingly orders $x^2$.
\end{proposition}
\begin{proof}
    Let us deduce the inequality $x^1_{\sigma^1(1)}-x^1_j\leq x^2_{\sigma^2(1)}-x^2_j$. The other one is similar. First, note  that we can use Proposition \ref{prop:MCMC_xvalues} and consider the following cases:
    \begin{itemize}
        \item $x^1_j=o_j$. This means that $x^2_j=o_j$ and consequently 
        \begin{gather*}
            x^1_{\sigma^1(1)}-x^1_j\leq x^2_{\sigma^2(1)}-x^2_j.
        \end{gather*}
        \item $x^1_j=x^1_{\sigma^1(1)}$. Then
        \begin{gather*}
            0=x^1_{\sigma^1(1)}-x^1_j\leq x^2_{\sigma^2(1)}-x^2_j
        \end{gather*}
        \item $x^1_j=x^1_{\sigma^1(n)}$. Now, if we apply again Proposition $\ref{prop:MCMC_xvalues}$, we obtain three subcases
        \begin{itemize}
            \item $x^2_j=o_j$. Then $x^1_{\sigma^1(n)}\geq o_j$ and 
            \begin{gather*}
                x^1_{\sigma^1(1)}-x^1_j=x^1_{\sigma^1(1)}-x^1_{\sigma^1(n)}\leq x^2_{\sigma^2(1)}-o_j=x^2_{\sigma^2(1)}-x^2_j
            \end{gather*}
            \item $x^2_j=x^2_{\sigma^2(n)}$. Then
            \begin{gather*}
                x^1_{\sigma^1(1)}-x^1_j=x^1_{\sigma^1(1)}-x^1_{\sigma^1(n)}=\delta_1<\delta_2=x^2_{\sigma^2(1)}-x^2_{\sigma^2(n)}=x^2_{\sigma^2(1)}-x^2_j
            \end{gather*}
            \item $x^2_j=x^2_{\sigma^2(1)}$. This case is unfeasible since $x^2_j=x^2_{\sigma^2(1)}$ implies $o_j\geq x^2_{\sigma^2(1)}$ and thus $|o_j-x_j^1|\geq |o_j-x^1_{\sigma(1)}|$, from which taking $x^1_j=x^1_{\sigma(1)}$ would produce a feasible point in $\mathcal{R}_{\delta_1}$ with lesser cost.
        \end{itemize}
    \end{itemize}
\end{proof}
The following result contributes to the satisfaction of the axioms in the previous proposition.
\begin{proposition}\label{P:nestedexist}
    Let $\delta_1,\delta_2\in[0,1]$ such that $\delta_1<\delta_2$, $c\in\io^n$ and $o\in\io^n$ verifying $\max\{o_1,...,o_n\}-\min\{o_1,...,o_n\}>\delta_2$. Consider $x^2\in\operatorname{arg}\min_{y\in\mathcal{R}_{\delta_2}}\xi_c^o(y)$. Then, we can find  $x^1\in\operatorname{arg}\min_{y\in\mathcal{R}_{\delta_1}}\xi_c^o(y)$ such that 
    \begin{gather*}
        x^2_{\sigma^2(n)}\leq x^1_{\sigma^1(n)}\leq  x^1_{\sigma^1(1)}\leq x^2_{\sigma^2(1)}.
    \end{gather*}
\end{proposition}
\begin{proof}
    We proceed by employing a contrapositive argument. Let us assume that we cannot find $x^1$ between the lower and upper bounds of $x^2$. This implies that either we have  $x^1_{\sigma^1(1)}>x^2_{\sigma^2(1)}$ or $x^1_{\sigma^1(n)}<x^2_{\sigma^2(n)}$. We will analyze the case $x^1_{\sigma^1(1)}>x^2_{\sigma^2(1)}$ and show that $x^1$ is not a solution. The proof of the case $x^1_{\sigma^1(1)}<x^2_{\sigma^2(1)}$ is analogous.
    
    Since $x^1_{\sigma^1(1)}>x^2_{\sigma^2(1)}$, $x^1\in\mathcal{R}_{\delta_1}$, $x^2\in\mathcal{R}_{\delta_2}$, then the maximum value $o^+=\max\{o_1,...,o_n\}$ must satisfy $ x^1_{\sigma^1(1)}>x^2_{\sigma^2(1)}\geq o^+$. However, if we consider
    \begin{gather*}                     x^0\in\operatorname{arg}\min_{\substack{y\in\mathcal{R}_{\delta_1}\\
    y_j\in[x^2_{\sigma^2(n)},x^2_{\sigma^2(1)}]}}\xi_c^o(y)
    \end{gather*}
    $x^0_{\sigma_0(1)}\leq x^2_{\sigma^2(1)}$ and thus $o^+\geq x^2_{\sigma^2(1)}$. Consequently, $o^+= x^2_{\sigma^2(1)}$.
    Now, think about the point $\hat{x}\in\io^n$ defined by
    \begin{gather*}
        \hat{x}_j=\begin{cases}
            x^1_j&\text{ if }x^1_j\neq x^1_{\sigma^1(1)}\\
            o^+&\text{ if }x^1_j = x^1_{\sigma^1(1)}
        \end{cases}.
    \end{gather*}
    Since $x^1_{\sigma^1(1)}>o^+$, $\hat{x}\in\mathcal{R}_{\delta_1}$. In addition, $\xi_c^o(\hat{x})<\xi_c^o(x^1)$ and then $x^1$ is not a solution to $\min_{y\in\mathcal{R}_{\delta_1}}\xi_c^o(y)$, which is a contradiction.
\end{proof}
The two previous results should be interpreted as a study on \textit{nested} solutions to MCMC, i.e., how decreasing the value of $\delta$ in the MCMC model modifies the values of the extreme opinions towards a more central vision. 

\section{Some considerations about OWA-MCC}\label{sec:OWAMCC}
In the previous section, we presented the mutual consensus measure. We have also proven that if the opinions of the members of a group satisfy the restriction defined by the mutual consensus, such opinions necessarily satisfy the constraints determined by other classical consensus measures. Later, we will use this fact in the context of the OWA-based MCC model \citep{zhang2011}, to gain some understanding of the relationship between this model and the MCMC model. Having this goal in mind, this section provides some considerations regarding the OWA-MCC model, sets notations, and gives some geometric properties of OWA-MCC. 
We focus on the most basic case in which the feasible region is determined by the consensus measure $\kappa^{\Psi_\omega}$, where $\Psi_\omega$ is the OWA operator with the weighting vector $\omega$.

\begin{definition}[OWA-MCC model] Let us consider $n$ DMs who provide their opinions $o=(o_1,...,o_n)\in\io^n$ using a numerical scale. Assume that $c=(c_1,...,c_n)\in\io^n$ is the cost vector indicating costs of modifying DMs' opinions, $\omega=(\omega_1,...,\omega_n)$ are the weights for the OWA operator, and $\varepsilon\in\io$ is the consensus threshold. Then, the OWA-MCC model can be formulated as follows:

    \begin{equation}\label{OWA-MCC}
    \tag{OWA-MCC}    \min_{x\in\mathcal{R}_\varepsilon^{\Psi_\omega}}\xi_c^o(x),
    \end{equation}
    where $\rowae=\{x\in\io^n\tq \kappa^{\Psi_\omega}(x)\leq \varepsilon\}$.
\end{definition}
The following result, in line with the study in \cite{paramcons}, provides some geometric properties of the feasible region of OWA-MCC under different consensus measures.
\begin{theorem}\label{sym-shift}
    Let us consider $\delta,\varepsilon,\gamma\in\io$, a weighting vector $\omega$ for the OWA operator, and another weighting vector $w$. Then, the following statements hold:
    \begin{enumerate}
        \item $\rmd$ is convex, symmetric, and shift-invariant.
        
        \item $\rowae$ is symmetric and shift-invariant.
        \item $\mathcal{R}^{w,\Psi_\omega}_\gamma$ is  shift-invariant. If  $w_1=w_2=...=w_n=\frac{1}{n}$, it is symmetric.
        \item $\mathcal{R}^{w}_\gamma$ is convex  and shift-invariant. If $w_1=w_2=...=w_n=\frac{1}{n}$, it is also symmetric.        
    \end{enumerate}
\end{theorem}
\begin{proof} 
    Below, we provide the proof for the second statement. The remaining ones are similar. 

    Let $x\in\rowae$ and consider a permutation $\sigma:\{1,2,...,n\}\to \{1,2,...,n\}$. It is clear that, for $y=x_{\sigma}$, $\Psi_\omega(x)=\Psi_\omega(y)$. Therefore, $|y_i-\Psi_\omega(y)|=|x_{\sigma^{-1}(i)}-\Psi_\omega(x)|\leq\varepsilon$ for $i=1,2,...,n$ and thus $y\in \rowae$. Additionally, if  $x\in\rowae$, and we pick $\lambda\in\R$ such that $y=(x_1+\lambda,x_2+\lambda,...,x_n+\lambda)\in\io^n$, then $\Psi_\omega(y)=\lambda +\Psi_\omega(x)\in\io$ and for any $i=1,2,...,n$, we obtain $|y_i-\Psi_\omega(y)|=|x_i-\Psi_\omega(x)| \leq\varepsilon$, which implies $y\in\rowae$.
\end{proof}
The knowledge about the geometry of OWA-MCC allows us to provide a general method for linearizing the OWA-MCC model under the different consensus measures, provided that a certain symmetry condition is met on the parameters.
\begin{theorem}\label{T:LOWA-Eq_cost_w}
    Let us consider a group of $n$ DMs who express their opinions $o = (o_1, \ldots, o_n) \in \io^n$ on a numerical scale. Assume that the cost of modifying each DM's opinion is uniform, i.e., $c_1 = c_2 = \cdots = c_n = \frac{1}{n}$. Let $\omega = (\omega_1, \ldots, \omega_n)$ denote the weight vector used in the OWA operator, and assume that the importance weights for the DMs are equal, i.e., $w_1 = w_2 = \cdots = w_n = \frac{1}{n}$.
Let us fix  the consensus thresholds $\varepsilon,\delta,\gamma_1,\gamma_2\in\io$, consider a permutation $\sigma:\{1,...,n\}\to\{1,...,n\}$ such that $o_{\sigma(1)}\geq o_{\sigma(2)}\geq...\geq o_{\sigma(n)}$ and define $\psi_\omega:\io^n\to\io$ as $\psi_\omega(x)=\sum_{k=1}^n\omega_kx_k \pt x\in\io^n$. If $x^*\in\io$ solves the linear programming problem
    \begin{equation*}
        \begin{split}
            &\min_{x\in\io^n}\xi_c^{o_\sigma}(x)\\
            \text{s.t.}&\begin{cases}         x\in\mathcal{R}_\delta\cap\mathcal{R}_\varepsilon^{\psi_\omega}\cap\mathcal{R}_{\gamma_1}^{w,\psi_\omega}\cap\mathcal{R}_{\gamma_2}^{w}\\
            x_1\geq x_2\geq...\geq x_n
            \end{cases}
        \end{split},
    \end{equation*}
    then, $x^*_{\sigma^{-1}}$ is a solution for the optimization problem
    \begin{equation*}
        \begin{split}
            &\min_{x\in\io^n}\xi_c^{o}(x)\\
            \text{s.t.}&\begin{cases}         x\in\mathcal{R}_\delta\cap\mathcal{R}_\varepsilon^{\Psi_\omega}\cap\mathcal{R}_{\gamma_1}^{w,\Psi_\omega}\cap\mathcal{R}_{\gamma_2}^{w}\\
            \end{cases}
        \end{split}.
    \end{equation*}
\end{theorem}
\begin{proof}
    If $x^*$ solves the first optimization problem, the previous results about symmetry combined with the fact that $\Psi_\omega(x)=\psi_\omega(x)\pt x\in\io$ such that $x_1\geq x_2\geq...\geq x_n$, lead to $x^*_{\sigma^{-1}}\in\mathcal{R}_\delta\cap\mathcal{R}_\varepsilon^{\Psi_\omega}\cap\mathcal{R}_{\gamma_1}^{w,\Psi_\omega}\cap\mathcal{R}_{\gamma_2}^{w}$.
    Now, consider a solution $y^*$ for the second optimization problem. Due to the symmetry of the feasible region, $y^*_\sigma$ belongs to the feasible region of the first optimization model. Since $x^*$ is the solution to the first problem, $\xi_c^o(y^*)=\xi_c^{o_\sigma}(y^*_\sigma)\geq\xi_c^{o_\sigma}(x^*)$. However, $y^*$ is the solution to the second one and thus $\xi_c^{o_\sigma}(x^*)=\xi_c^{o}(x^{*}_{\sigma^{-1}})\geq \xi_c^o(y^*)$, 
    which implies that $\xi_c^o(y^*)=\xi_c^{o_\sigma}(x^*)$. 
\end{proof}

This result generalizes the linearization provided by \cite{zhang2011} and provides a linear programming formulation that allows solving the OWA-MCC model, which is a nonlinear optimization problem in the general case. However, it strongly depends on the symmetry conditions $c_1=c_2=...=c_n=\frac{1}{n}$ and $w_1=w_2=,...=w_n=\frac{1}{n}$. When we drop these assumptions, it is hard to obtain a linear reformulation for the general MCC problem under OWA aggregations.

Let us consider $n$ DMs who provide their opinions $o=(o_1,...,o_n)\in\io^n$. Assume that $c=(c_1,c_2,...,c_n)\in \io^n$ are the costs of modifying their opinions, $\omega=(\omega_1,...,\omega_n)\in \io^n$ are the weights for the OWA operator, and  $w=(w_1,w_2,...,w_n)\in \io^n$ are the weights for the consensus meausures. If we fix the consensus thresholds $\varepsilon,\delta,\gamma_1,\gamma_2\in\io$, we obtain the general model 
\begin{equation*}
\tag{G-OWA-MCC}
        \begin{split}
            &\min_{x\in\io^n}\xi_c^{o}(x)\\
            \text{s.t.}&\begin{cases}         x\in\mathcal{R}_\delta\cap\mathcal{R}_\varepsilon^{\Psi_\omega}\cap\mathcal{R}_{\gamma_1}^{w,\Psi_\omega}\cap\mathcal{R}_{\gamma_2}^{w}\\
            \end{cases}
        \end{split}.
    \end{equation*}

\begin{figure}[H]
        \centering
        \begin{subfigure}[b]{0.45\textwidth}  
             \centering
            \includegraphics[width=\textwidth]{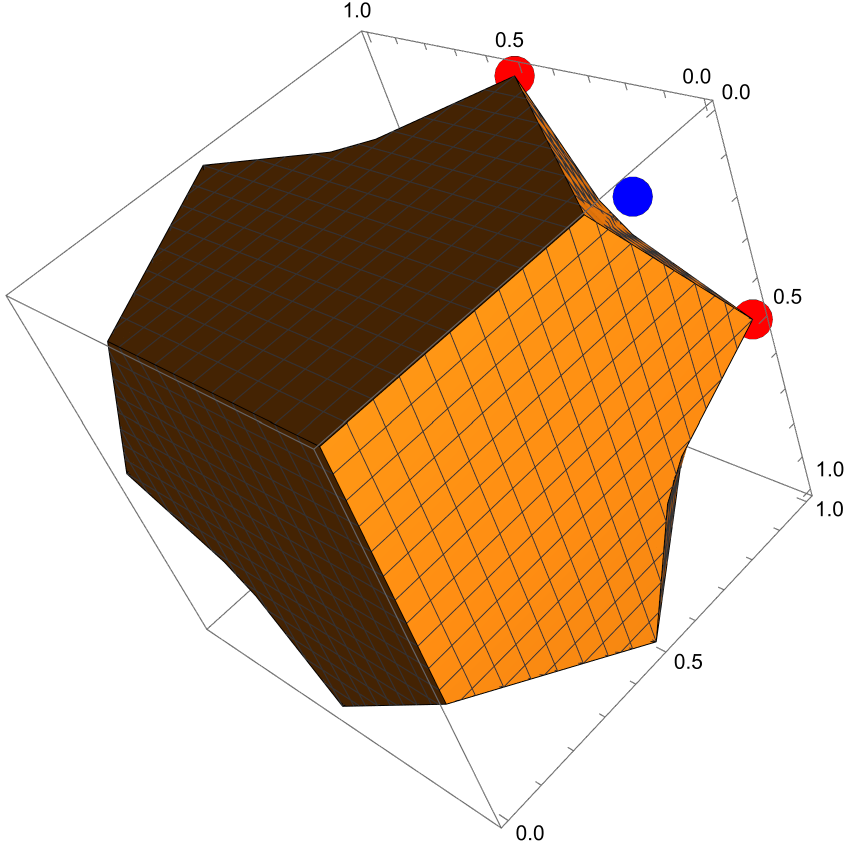}            \caption{Draft of the region $\mathcal{R}_\frac{1}{2}^{\Psi_\omega}$ for $\omega=(0.25,0.5,0.25)$.}
            \label{img:epsnonconv}
        \end{subfigure}\hfill
        \begin{subfigure}[b]{0.45\textwidth}
        \centering
        \includegraphics[width=0.9\textwidth]{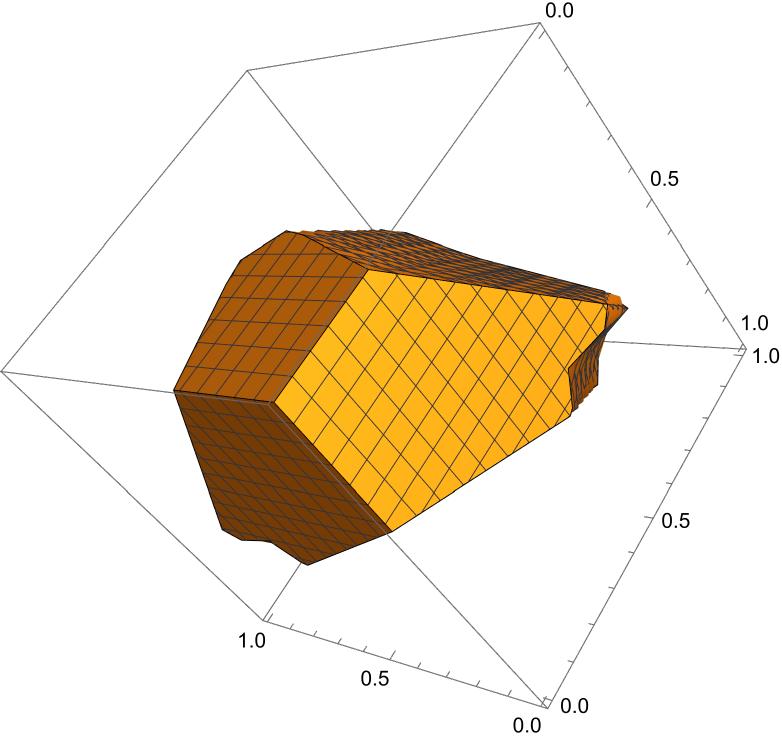}
        \caption{Draft of the region $\mathcal{R}_\frac{1}{2}^{w,\Psi_\omega}$ for $\omega=(0.7,0.1,0.2)$ and $w=(0.25,0.25,0.5)$.}
            \label{img:gammanonconv}
        \end{subfigure}
        \caption{Examples showing the non-convexity of the regions that determine the G-OWA-MCC model.}
    \end{figure}

    Note that the absence of the aforementioned symmetry conditions implies that the region $\rowae$ is not convex in general.
    For instance, let us fix $n=3$,  $\varepsilon=\frac{1}{2}$ and $\omega=(\frac{1}{4},\frac{1}{2},\frac{1}{4})$. The points $x_1=(1,\frac{1}{2},0)$ and $x_2=(1,0,\frac{1}{2})$ satisfy $\Psi_\omega(x_1)=\Psi_\omega(x_2)=\frac{1}{2}$ and thus $x_1,x_2\in \mathcal{R}_\frac{1}{2}^{\Psi_\omega}$. However, the middle point $m=\frac{1}{2}(x_1+x_2)=(1,\frac{1}{4},\frac{1}{4})$, satisfies $\Psi_\omega(m)=\frac{7}{16}$ and thus $|1-\Psi_\omega(m)|=\frac{9}{16}\geq \frac{1}{2}$, which means that $m\notin\mathcal{R}_\frac{1}{2}^{\Psi_\omega}$ (see Fig. \ref{img:epsnonconv}). Furthermore, for the configuration $\omega=(0.7,0.1,0.2), w=(0.25,0.25,0.5)$ and $\gamma=0.2$, it is easy to check that the region $\mathcal{R}_\gamma^{w,\Psi_\omega}$ is non-convex and non-symmetric (see Fig. \ref{img:gammanonconv}). This non-convexity implies that the feasible region is not a convex polytope \citep{paramcons}, which suggests that the model is inherently non-linear. This makes the resolution of the model G-OWA-MCC extremely costly in LSGDM scenarios involving hundreds or thousands of DMs.

\section{Relationship between Mutual Consensus and OWA-based Consensus measures}\label{sec:relationship}
In the previous section, we showed that some symmetry conditions allow the linearization of the OWA-MCC model. However, in the general case, the feasible region in the OWA-MCC model may consist of a non-convex polytope, which complicates the computational resolution of the optimization problem. A classical approach for solving it consists of the use of binary linear programming \citep{zhang2013}, but this vision may not be applied to solve LSGDM problems due to the elevated computational time required. We aim to apply the linear MCMC model to provide an approximate solution for the general OWA-MCC problem. In this section, we study the relationship between the models OWA-MCC and MCMC.

Let us start by analyzing the conditions under which, given $\varepsilon\in\io$, we can find $\delta_-,\delta_+\in\io$ such that $\cR_{\delta_-}\subseteq\cR_{\varepsilon}^{\Psi_\omega}\subseteq\cR_{\delta_+}$. 
\begin{theorem} Let us consider the weighting vector $\omega=(\omega_1,...,\omega_n)\in\io^n$. Then $\rmd\subseteq\rowae$ if and only if $\delta\leq \min\{\frac{\varepsilon}{1-\min\{\omega_1,\omega_n\}},1\}$.    
\end{theorem}
\begin{proof}
    We omit the case $\frac{\varepsilon}{1-\min\{\omega_1,\omega_n\}}\geq1$. Assume that $\rmd\subseteq\rowae$. In particular, the vertices of $\rmd$ in the form $x\in\{0,\delta\}^n$ should also be elements of $\rowae$. Let $\vartheta=1,2,...,n-1$ be the number of zeros in the components of $x$. Note that we intentionally exclude the cases $\vartheta\in\{0,n\}$ since they provide no information. Then $\Psi_\omega(x)=\delta\sum_{k=\vartheta+1}^nw_k$. Using that $x\in\rowae$, we obtain that
    \begin{equation*}
        \begin{split}
            \varepsilon&\geq|\delta-\delta\sum_{k=\vartheta+1}^n\omega_k|=\delta\sum_{k=1}^\vartheta w_k\\
            \varepsilon&\geq|0-\delta\sum_{k=\vartheta+1}^n\omega_k|=\delta \sum_{k=\vartheta+1}^n\omega_k
        \end{split}
    \end{equation*}
    and then $\delta\max\{\sum_{k=1}^\vartheta w_k,\sum_{k=\vartheta+1}^nw_k\}\leq\varepsilon$ for any value $\vartheta=1,2,...,n-1$. This, implies that $\delta\max\{\max_{\vartheta=1,...,n}\{\sum_{k=1}^\vartheta \omega_k\},\max_{\delta=1,...,n}\sum_{k=\vartheta+1}^n\omega_k\}\}\leq \varepsilon$, i.e., $\delta\max\{1-\omega_1,1-\omega_n\}\leq \varepsilon$.

    Conversely, assume that $\delta\leq \frac{\varepsilon}{1-\min\{\omega_1,\omega_n\}}$ and let us analyze the points of the form $x=(\delta, y,0)$, with $y\in [0,\delta]^{n-2}$. Note that $\Psi_\omega(x)=\delta\omega_1+\sum_{k=2}^{n-1}\omega_ky_k$  and 
    \begin{equation*}
    \begin{split}
        |x_1-\Psi_\omega(x)|&=\delta-\delta\omega_1-\sum_{k=2}^{n-1}\omega_ky_k=\delta\sum_{k=1}^{n-1}\omega_k-\sum_{k=1}^{n-1}\omega_kx_k=\\
        &=\sum_{k=1}^{n-1}\omega_k(\delta-x_k)\leq\delta(1-\omega_1)\\
        |\Psi_\omega(x)-x_n|&=\delta\omega_1+\sum_{k=2}^{n-1}\omega_ky_k\leq \delta \sum_{k=1}^{n-1}\omega_k=\delta(1-\omega_n)
    \end{split}
    \end{equation*}
    Now, due to $\delta\leq \frac{\varepsilon}{1-\min\{\omega_1,\omega_n\}}$, we obtain
    \begin{equation*}
        \kappa^{\Psi_{\omega}}(x)\leq\max\{|x_1-\Psi_\omega(x)|,|x_n-\Psi_\omega(x)|\}\leq\delta(1-\min\{\omega_1,\omega_n\})\leq\varepsilon.
    \end{equation*}
    Consequently, the points $x=(\delta,y,0)$, which characterize the boundary of $\rmd$ are contained in $\rowae$. Using the convexity and shift-invariance of $\rmd$, we conclude that $\rmd\subseteq\rowae$.
\end{proof}



\begin{figure}[H]
    \centering
    \includegraphics[width=0.57\linewidth]{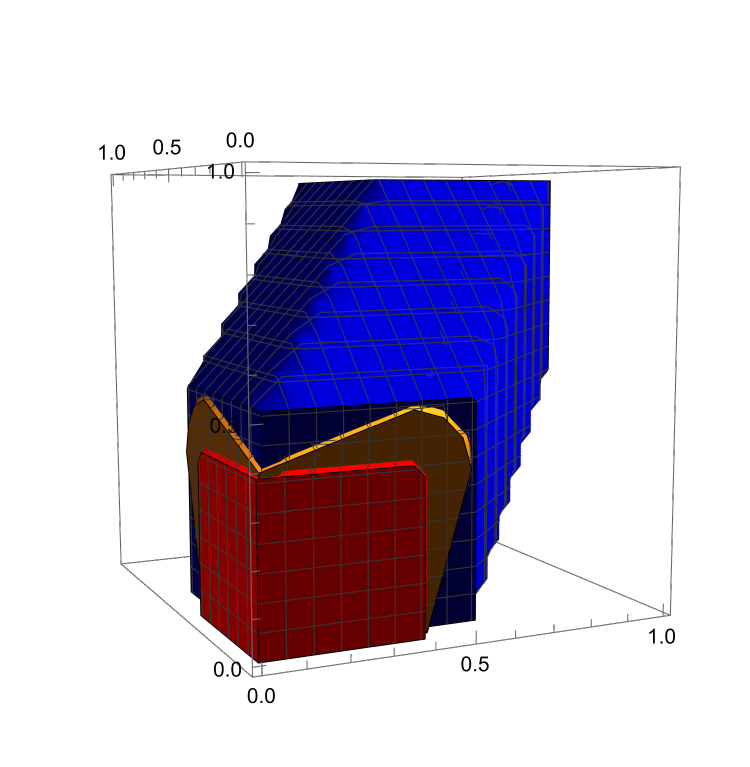}
    \caption{Draft of the regions $\mathcal{R}_{\delta_-},\; \rowae$ and $\mathcal{R}_{\delta_+}$ for $\varepsilon=0.25$ and $\omega=(0.35,0.2,0.45)$.}
    \label{fig:example-epsdelta}
\end{figure}

\begin{theorem} Let us consider a weighting vector $\omega=(\omega_1,...,\omega_n)\in\io^n$. Then, $\rowae\subseteq\rmd$ if and only if $\delta\geq\min\{1,2\varepsilon\}$.  
\end{theorem}
\begin{proof}
We will assume that $2\varepsilon<1$ because the case $\varepsilon\geq\frac{1}{2}$ can be easily obtained. 
Suppose first that $\delta\geq2\varepsilon$ and pick $x\in\mathcal{R}_\varepsilon^{\Psi_\omega}$ such that $x_1\geq x_2\geq...\geq x_n$. This, implies that $\max\{|x_1-\Psi_\omega(x)|, |\Psi_\omega(x)-x_n|\}\leq\varepsilon$. Consequently, 
\begin{equation*}
    |x_i-x_j|\leq|x_1-x_n|=| x_1-\Psi_\omega(x)|+|\Psi_\omega(x)-x_n|\leq 2\varepsilon\leq\delta.
\end{equation*}
Now, let us suppose that $\rowae\subseteq\rmd$, and define $k_0=\min\{k=1,...,n\tq \sum_{j=1}^{k-1} \omega_j<\frac{1}{2}, \sum_{j=1}^{k} \omega_j\geq\frac{1}{2}\}$. At this stage, consider the vector 
\begin{equation*}
    x=(2\varepsilon,2\varepsilon,....,2\varepsilon, \underbrace{\beta}_{k_0},0,...,0)
\end{equation*}
where $\beta=\varepsilon\frac{1-2 \sum_{j=1}^{k_0-1} \omega_j}{\omega_{k_0}}$. Note that the conditions over $k_0$ guarantee that $\beta\in[0,2\varepsilon]$. Consequently, $\Psi_\omega(x)=2\varepsilon\sum_{j=1}^{k_0-1} \omega_j+\omega_{k_0}\beta=\varepsilon$. This implies that $|x_k-\Psi_\omega(x)|\leq\varepsilon\pt k=1,...,n$ and then $x\in\rowae$. Since by hypothesis $x\in\rmd$, it is necessary that $|x_1-x_n|=2\varepsilon\leq\delta$.
\end{proof}
Let us recall the two previous results in a single Theorem.
\begin{theorem}\label{T:deltabounds}
    Let $\omega=(\omega_1,...,\omega_n)\in\io^n$ be a weighting vector. Then, for any $\varepsilon\in\io$, the value $ \delta_-=\min\{\frac{\varepsilon}{1-\min\{\omega_1,\omega_n\}},1\}$ satisfies 
    \begin{equation*}
        \delta_-=\max\{\delta\in\io\tq\rmd\subseteq\rowae\}, 
    \end{equation*}
    whereas, the value $\delta_+=\min\{2\varepsilon,1\}$ satisfies
    \begin{equation*}
        \delta_+=\min\{\delta\in\io\tq \rowae\subseteq\rmd\}.
    \end{equation*}    
\end{theorem}
Note that the previous theorem states that $\delta_-$ determines the biggest region determined by the mutual consensus measure $\kappa$ that is completely contained in $\rowae$, whereas $\delta_+$ defines the smallest region associated with $\kappa$ that contains $\rowae$ (see Figure \ref{fig:example-epsdelta}).
The most immediate consequence of this result provides a first approximation for a solution to OWA-MCC. We know that the solution to OWA-MCC, if not trivial, belongs to the region $\cR_{\delta_+}\setminus(\cR_{\delta_-})^\circ$, where $A^\circ$ denotes the topological interior of $A\subset\io$. Therefore, a solution $x^-\in\io^n$ to the MCMC problem $\min_{x\in\cR_{\delta_-}} \xi_c^o(x)$ satisfies the consensus constraint $\kappa^{\Psi_\omega}(x^-)\leq\varepsilon$ and provides the upper bound $\xi_c^o(x^-)$ for the cost of OWA-MCC, becoming an approximate solution for the worst-case scenario. On the other hand, a solution $x^+\in\min_{x\in\cR_{\delta_+}} \xi_c^o(x)$ provides a lower bound for the cost $\xi_c^o(x^+)$. Let us highlight this in a corollary.
\begin{corollary}\label{cor:costbound}
    Let $x^-\in\arg\min\limits_{x\in\cR_{\delta_-}} \xi_c^o(x)$ and $x^+\in\arg\min\limits_{x\in\cR_{\delta_+}} \xi_c^o(x)$. If  $x\in\arg\min\limits_{x\in\rowae} \xi_c^o(x)$, then $\xi_c^o(x)\in[\xi_c^o(x^+),\xi_c^o(x^-)]$.
\end{corollary}

Below, we show that the higher the number of variables/DMs involved, the lower the volumes of the regions $\cR_{\delta_-}$, $\rowae$, and $\cR_{\delta_+}$.  

\begin{corollary}\label{coro:conv}
    Let $\varepsilon\in]0,\frac{1}{2}[$, $\N^*:=\N\setminus\{1\}$ and consider a sequence of weighting vectors $\omega:\N^*\to\cup_{k=1}^\infty\R^k$ such that $\omega^n\in\io^n\; \forall n\in\N^*$. Let us denote
    \begin{equation*}
        \begin{split}
            \mathcal{R}^n_-&=\{x\in[0,1]^n\tq\kappa(x)\leq\frac{\varepsilon}{1-\min\{\omega_1^n,\omega_n^n\}}\}\\ \mathcal{R}^n_\varepsilon&=\{x\in[0,1]^n\tq\kappa^{\Psi_{\omega^n}}(x)\leq\varepsilon\}\\ \mathcal{R}^n_+&=\{x\in[0,1]^n\tq\kappa(x)\leq2\varepsilon\}
        \end{split}, \pt n\in\N^*.
    \end{equation*}
    
    Then, 
    \begin{equation*}
        \lim_{n\to\infty}\text{vol}(\mathcal{R}^n_-)=\lim_{n\to\infty}\text{vol}(\mathcal{R}^n_\varepsilon)=\lim_{n\to\infty}\text{vol}(\mathcal{R}^n_+)=0.
    \end{equation*}
\end{corollary}
\begin{proof}
    First, note that for $n\in\N^*$, $\min\{\omega_1^n,\omega_n^n\}\leq\frac{1}{n}\leq\frac{1}{2}$. Consequently, since $\varepsilon\leq\frac{1}{2}$, 
    \begin{equation*}
        \frac{\varepsilon}{1-\min\{\omega_1^n,\omega_n^n\}}\leq 1
    \end{equation*}
    and thus $\mathcal{R}^n_-, \mathcal{R}^n_\varepsilon$, and $\mathcal{R}^n_+$ stand respectively for the regions $\mathcal{R}_{\delta_-},\rowae$ and $\mathcal{R}_{\delta_+}\pt n\in\N^*$ mentioned in the previous result. Now, given $n\in\N^*$, it is clear that $\text{vol}(\mathcal{R}^n_+)=2^n\varepsilon^n$. Since $\mathcal{R}^n_-\subseteq \mathcal{R}^n_\varepsilon\subseteq\mathcal{R}^n_+\pt n\in\N^*$ and $\lim_{n\to\infty}2^n\varepsilon^n=0$, the conclusion follows.
\end{proof}
In particular $\lim_{n\to\infty}\text{vol}(\mathcal{R}^n_\varepsilon)-\lim_{n\to\infty}\text{vol}(\mathcal{R}^n_-)=0$, which is a kind of convergence result.
\section{An MCMC-based algorithm to approximate the solution of OWA-MCC}\label{sec:algorithm}
In the previous section, we found that $x^-\in\arg\min_{x\in\cR_{\delta_-}} \xi_c^o(x)$ satisfies the consensus constraint $\kappa^{\Psi_\omega}(x^-)\leq\varepsilon$, which makes it an approximate solution for OWA-MCC. However, the cost value for the OWA-MCC problem belongs to the interval $[\xi_c^o(x^+),\xi_c^o(x^-)]$. In this section, we aim to refine this result. To do so, given $\varepsilon\in]0,\frac{1}{2}[$, we will find a value $\delta_*$ such that $\cR_{\delta_-}\subseteq \cR_{\delta_*}\subseteq \cR_{\delta_+}$ and a solution $x^*\in\arg\min_{x\in\cR_{\delta_*}} \xi_c^o(x)$ that satisfies $\xi_c^o(x^*)\in[\xi_c^o(x^+),\xi_c^o(x^-)]$ and $\kappa^{\Psi_\omega}(x^*)=\varepsilon$. Below, we describe the computational process, as well as the formal results that sustain the computations.

Let us fix $\varepsilon\in]0,\frac{1}{2}[$ and compute the values $\delta_-$ and $\delta_+$ as in Theorem \ref{T:deltabounds}. Assume that $o\notin\rowae$. Since $\delta_-\leq \varepsilon$, we may use Proposition \ref{P:MC-Inside} and deduce that $\kappa^{\Psi_\omega}(x^-)\leq\varepsilon$, for $x^-\in\arg\min_{x\in\cR_{\delta_-}} \xi_c^o(x)$. Additionally, Proposition \ref{P:nestedexist} allow us taking $x^+\in\arg\min_{x\in\cR_{\delta_+}} \xi_c^o(x)$ such that $x^-$ and $x^+$ are nested. Before we advance, we need the following complementary result. 
\begin{proposition}
    Let $\varepsilon\in\io$ and $\alpha\geq2\varepsilon$. If $x\in\cR_\alpha$ and $x_{\sigma(1)}-x_{\sigma(n)}=\alpha$, were $\sigma$ is the permutation that decreasingly orders $x$, then $\kappa^\Psi(x)\geq\varepsilon$.
\end{proposition}
\begin{proof}
    Let us suppose that  $\kappa^\Psi(x)<\varepsilon$. Let us denote 
     \begin{gather*}         \beta=\sum_{i=1}^n\omega_i(x_{\sigma(1)}-x_{\sigma(i)})
     \end{gather*}
     to deduce that
     \begin{gather*}
         \kappa^\Psi(x)=\max\{\sum_{i=1}^n\omega_i(x_{\sigma(1)}-x_{\sigma(i)}),\sum_{i=1}^n\omega_i(x_{\sigma(i)}-x_{\sigma(n)})\}=\max\{\beta,\alpha-\beta\}
     \end{gather*}
     and consequently, since $\kappa^\Psi(x)<\varepsilon$, it must be $\max\{\beta,\alpha-\beta\}<\varepsilon$. However, this implies that $\beta<\varepsilon$ and $\alpha-\beta<\varepsilon$. The last inequality can be expressed as $\beta>\alpha-\varepsilon\geq\varepsilon$, which is a contradiction.
\end{proof}

If $o\notin\rowae$, then $x^+$ belongs to $\cR_{\delta_+}$ and satisfies $x^+_{\sigma(1)}-x^+_{\sigma(n)}=\delta_+$. Thus, the previous property ensures that
\begin{gather*}
    \kappa^{\Psi_\omega}(x^-)\leq\varepsilon\leq\kappa^{\Psi_\omega}(x^+).
\end{gather*}
Let us show that, in non-trivial examples, we have $\kappa^{\Psi_\omega}(x_+)-\kappa^{\Psi_\omega}(x_-)\neq 0$.
\begin{proposition}
    Let us consider $\varepsilon\in]0,\frac{1}{2}[$, a cost vector $c$, an initial opinion $o$, and a weighting vector $\omega$, and define $\delta_-=\min\{\frac{\varepsilon}{1-\min\{\omega_1,\omega_n\}},1\}$ and $\delta_+=2\varepsilon$. If $x_-\in\arg\min_{x\in\mathcal{R}_{\delta_-}}\xi_c^o(x)$ and $x_+\in\arg\min_{x\in\mathcal{R}_{\delta_+}}\xi_c^o(x)$ are such that $\kappa^{\Psi_\omega}(x_+)=\kappa^{\Psi_\omega}(x_-)$, then $x_+\in\arg\min_{x\in\mathcal{R}_\varepsilon^{\Psi_\omega}}\xi_c^o(x)$.
\end{proposition}
\begin{proof}
    Since $x_-\in\mathcal{R}_{\delta_-}$, $\kappa^{\Psi_\omega}(x_-)\leq\varepsilon$. Therefore, $\kappa^{\Psi_\omega}(x_+)\leq\varepsilon$, and $x_+\in\mathcal{R}_\varepsilon^{\Psi_\omega}$.
    Since $x_+\in\arg\min_{x\in\mathcal{R}_{\delta_+}}\xi_c^o(x)$, and $\mathcal{R}_\varepsilon^{\Psi_\omega}\subseteq\mathcal{R}_{\delta_+}$ then $x_+\in\arg\min_{x\in\mathcal{R}_\varepsilon^{\Psi_\omega}}\xi_c^o(x)$.
\end{proof}

Therefore, we have $\kappa^{\Psi_\omega}(x^-)\leq\varepsilon\leq\kappa^{\Psi_\omega}(x^+)$ with $\kappa^{\Psi_\omega}(x^-)<\kappa^{\Psi_\omega}(x^+)$. Let us show that, indeed, this happens whenever we have nested solutions $x^1$ and $x^2$ for two different values $\delta_1$ and $\delta_2$, respectively. 

\begin{theorem}
    Let $\delta_1,\delta\in[0,1]$ such that $\delta_1<\delta_2$, a weighting n-tuple $\omega\in\io^n$, $c\in\io^n$ and $o\in\io^n$ verifying $\max\{o_1,...,o_n\}-\min\{o_1,...,o_n\}>\delta_2$. Consider $x^1\in\operatorname{arg}\min_{y\in\mathcal{R}_{\delta_1}}\xi_c^o(y)$ and $x^2\in\operatorname{arg}\min_{y\in\mathcal{R}_{\delta_2}}\xi_c^o(y)$ such that 
    \begin{gather*}
         x^2_{\sigma^2(n)}\leq x^1_{\sigma^1(n)}\leq  x^1_{\sigma^1(1)}\leq x^2_{\sigma^2(1)}.
    \end{gather*}Then, $\kappa^{\Psi_\omega}(x^1)\leq \kappa^{\Psi_\omega}(x^2)$.
\end{theorem}
\begin{proof}
    First, note that given $x\in\io^n$, we can write
    \begin{gather*}
        \kappa^{\Psi_\omega}(x)=\max_{j=1,...,n}\{|x_j-\sum_{i=1}^n\omega_ix_{\sigma(i)}|\}=\max_{j=1,...,n}\{|\sum_{i=1}^n\omega_ix_j-\sum_{i=1}^n\omega_ix_{\sigma(i)}|\}=\\\max_{j=1,...,n}\{|\sum_{i=1}^n\omega_i(x_j-x_{\sigma(i)})|\}=\max\{\sum_{i=1}^n\omega_i(x_{\sigma(1)}-x_{\sigma(i)}),\sum_{i=1}^n\omega_i(x_{\sigma(i)}-x_{\sigma(n)})\}.
    \end{gather*}
    Now, the conclusion follows from Proposition \ref{P:d1d2MCM}.
\end{proof}

Let us recall the situation. We have $\kappa^{\Psi_\omega}(x^-)\leq\varepsilon\leq\kappa^{\Psi_\omega}(x^+)$ with $\kappa^{\Psi_\omega}(x^-)<\kappa^{\Psi_\omega}(x^+)$ and, whenever another pair $\delta_1$, $\delta_2$ generates nested solutions, we obtain that the consensus measure $\kappa^{\Psi_\omega}$ becomes monotonic. Since $\kappa^{\Psi_\omega}$ is also continuous, we can conjecture that for some value $\delta_*\in[\delta_-,\delta_+]$, we can find a solution $x^*\in\arg\min_{x\in\cR_{\delta_*}} \xi_c^o(x)$ such that $\kappa^{\Psi_\omega}(x^*)\leq\varepsilon$. In such a case, $x^*$ also solves the problem 
\begin{gather*}
    \min_{x\in\cR_{\delta_*}\cap\rowae} \xi_c^o(x)
\end{gather*}
and thus it will be the approximation we are looking for. So, finding an approximate solution to OWA-MCC is equivalent to computing the value $\delta_*$. To do so, we use the almost linear geometry of the involved constraints and functions to presume that the faster way to find such a $\delta^*$ could be through a proportional numeric method. So, we will use an iterative algorithm consisting of computing a provisional 
\begin{gather*}
    \delta_*=\frac{\varepsilon-\kappa^{\Psi_\omega}(x_-)}{\kappa^{\Psi_\omega}(x_+)-\kappa^{\Psi_\omega}(x_-)}(\delta_+-\delta_-)+\delta_-.
\end{gather*}
Now, we check the value $\kappa^{\Psi_\omega}(x^*)$ for the provisional value $x^*\in\arg\min_{x\in\cR_{\delta_*}}\xi_c^o(x)$, and compare $\kappa^{\Psi_\omega}(x^*)$ with both $\kappa^{\Psi_\omega}(x^-)$ and $\kappa^{\Psi_\omega}(x^+)$. In the following iteration, we may repeat the process by substituting either $\delta_-$ or $\delta_+$ with the recently computed $\delta_*$. We can repeat the process until $\kappa^{\Psi_\omega}(x^*)$ is close enough to $\varepsilon$. For the sake of clarity, we recall the main steps of this iterative method, which we call {\sc{ApOWAMCC}}, in Algorithm \ref{alg:approx}.

\begin{algorithm}[H]
\caption{Method to find an approximate solution to the OWA-MCC model}
\begin{algorithmic}[1]
\Procedure{ApOWAMCC}{$c,o,\omega,\varepsilon,N,\tau$}\Comment{Approximate the solution}
\State $\delta_-\gets \min\{\frac{\varepsilon}{1-\min\{\omega_1,\omega_n\}},1\}$
\State $x_-\gets \text{\sc{SolveMCMC}}(c,o,\delta_-)$
\State $\delta_+\gets \min\{2\varepsilon,1\}$
\State $x_+\gets \text{\sc{SolveMCMC}}(c,o,\delta_+)$
\State $n\gets 1$ 
\While{$(n\leq N) \& (\lvert \kappa^{\Psi_\omega}(x_-)-\varepsilon  \lvert>\tau)$}\Comment{We use $N$ as maximum number of iterations and $\tau$ as a tolerance}
\If{$\kappa^{\Psi_\omega}(x_+)=\kappa^{\Psi_\omega}(x_-)$}
\State $x_-\gets x_+$
\State $\delta_-\gets \delta_+$
\State $n\gets N+1$
\Else
\State $\delta\gets \frac{\varepsilon-\kappa^{\Psi_\omega}(x_-)}{\kappa^{\Psi_\omega}(x_+)-\kappa^{\Psi_\omega}(x_-)}(\delta_+-\delta_-)+\delta_-$
\State $x\gets \text{\sc{SolveMCMC}}(c,o,\delta)$
\If{$\kappa^{\Psi_\omega}(x)\leq\varepsilon$}
    \State $\delta_-\gets \delta$
    \State $x_-\gets x$
    \Else
    \State $\delta_+\gets \delta$
    \State $x_+\gets x$
\EndIf
    \State $n\gets n+1$
\EndIf
\EndWhile
\State \textbf{return} $x_-,\delta_-$\Comment{The approximate solution is $x_-$}
\EndProcedure
\end{algorithmic}\label{alg:approx}
\end{algorithm}
\section{Examples and Simulations}\label{sec:comput}
This section aims to illustrate the applicability of the theoretical tools developed in the manuscript. All the experiments were conducted using JuMP (Julia for Mathematical Programming) \citep{Dunning_2017}, a modeling language for mathematical optimization embedded in Julia \citep{Julia-2017}. The experiments were run on Julia 1.11.2 using Clp for linear programming and GLPK for BLP as optimizers, on a laptop with Windows 11 Home, a 2.2 GHz Intel Core i7-1360P CPU, and 16 GB of RAM. 

First, let us show two examples to illustrate the performance of \textsc{ApOWAMCC}. Each expert expresses an individual opinion, which is to be aggregated into a consensus value that minimizes the total cost under a constraint of the form $\kappa^{\Psi_\omega}\leq \varepsilon$. In each example, we compare our method with the BLP formulation \citep{zhang2013}. 
\begin{example}
    Let us consider a scenario involving five decision makers. The cost vector associated with modifying each expert’s opinion is given by ${c} = (1,\ 4,\ 3,\ 5,\ 2)/15$, indicating the relative difficulty of modifying each opinion. The aggregation is performed using an OWA operator with weights $\omega = (0.375,\ 0.1875,\ 0.25,\ 0.0625,\ 0.125)$. The original individual opinions are normalized to the unit interval and given by ${o} = (0.05,\ 0.1,\ 0.25,\ 0.3,\ 0.6)$. The maximum allowable deviation from the group opinion is set to $\varepsilon= 0.2$.
\end{example}
    If we apply \textsc{ApOWAMCC}, the resulting consensus solution is ${o^1} = (0.1,\ 0.1,\ 0.25,\allowbreak\ 0.3,\ 0.4333)$. The group consensus value derived is $0.3$, and the individual absolute deviations from this value are $(0.2,\ 0.2,\ 0.05,\ 0.0,\ 0.1333)$, all within the $\varepsilon$-bound. The total cost incurred to reach this solution is $0.0256$. This value lies within the cost range of $[0.01667,\ 0.03952]$ given by Corollary \ref{cor:costbound}. The computational time required for this method was $3$ milliseconds.

    The second method \citep{zhang2013}, applied to the same input data, yielded an identical consensus solution: ${o^2} = (0.1,\ 0.1,\ 0.25,\ 0.3,\ 0.4333)$. The group consensus value, individual deviations, and total cost remained the same as in the first method. However, the execution time was notably different, with this method requiring $16$ milliseconds to complete the computation.

    Although both methods arrived at the same consensus solution and incurred the same cost, the comparison highlights a difference in computational efficiency. The first method achieved the result in significantly less time, suggesting it may be more suitable for real-time or large-scale applications where computational speed is a critical factor.

\begin{example}
    This example examines the performance of the two previous methods for solving a GDM problem involving eight experts. The input opinions are normalized as \({o} = (0.05,\ 0.1,\ 0.25,\allowbreak\ 0.3,\ 0.6,\ 0.7,\ 0.5,\ 0.8)\), and the maximum allowable deviation from the group consensus is fixed at \(\varepsilon = 0.1\). The cost associated with modifying each opinion is given by ${c} = (0.0323,\ 0.1290,\allowbreak\ 0.0968,\ 0.1613,\ 0.0645,\ 0.1935,\ 0.0323,\ 0.2903)$. The aggregation of opinions is guided by an OWA operator with weights \(\omega = (0.175,\ 0.2,\ 0.0875,\ 0.25,\ 0.0325,\ 0.125,\allowbreak\ 0.1,\ 0.03)\).
\end{example}
    
    In the first method, the consensus opinion obtained is \(o^1=(0.5302,\ 0.5302,\ 0.5302,\ 0.5302,\allowbreak\ 0.6,\ 0.7,\ 0.5302,\ 0.7)\), leading to a group consensus value of \(0.6\). The individual absolute deviations from the consensus are \((0.0698,\ 0.0698,\ 0.0698,\ 0.0698,\ 0.0,\ 0.1,\ 0.0698,\ 0.1)\), all of which satisfy the \(\varepsilon\)-constraint. The total cost of achieving this solution is \(0.1653\), and it falls within the cost range \([0.1516,\ 0.1954]\) obtained by Corollary \ref{cor:costbound}. The computation was completed with a recorded execution time of \(0\) milliseconds.

    The second method \citep{zhang2013} generated a different consensus solution: \(o^2=(0.5,\ 0.5,\ 0.5,\ 0.5,\ 0.6,\ \allowbreak 0.7,\ 0.565,\ 0.7)\). The group consensus value remains at \(0.6\), but the individual deviations are now \((0.1,\ 0.1,\ 0.1,\ 0.1,\ 0.0,\ 0.1,\ 0.035,\ 0.1)\). This solution complies with the \(\varepsilon\)-bound for all the experts, and the total cost is lower, amounting to \(0.1537\), suggesting a more efficient reallocation of opinions in terms of cost. However, this result required significantly more time to compute, with an execution time of \(1432\) milliseconds.

    We can see that the second method achieved a lower consensus cost, but at the expense of substantially longer computation time. Conversely, the first method provided a near-optimal solution extremely quickly, which suggests that it may be preferable in applications where time efficiency is critical.

    Given these results, subsequently we aim to compare the solution provided by our algorithm with the exact solution provided under the symmetric setting described in Theorem \ref{T:LOWA-Eq_cost_w} and the classical BLP-based formulation \citep{zhang2013}. We measure the error in the objective function and the computational time (in milliseconds) required in each resolution method. In each test, we carry out $100$ random simulations and then we compute the mean and standard deviation of these simulations for different numbers of DMs.

The parameters are set as follows:
\begin{itemize}
    \item The initial opinions $o$ and the OWA weights $\omega$ are randomly generated. The consensus parameter is $\varepsilon=0.15$.
    \item The cost vector $c$ is fixed and symmetric in the Linear version, and random in the BLP version. In both cases, it has been normalized by dividing by the sum of its components.
    \item The maximum number of iterations in the numeric method to estimate $\delta_*$ is $N=10$, whereas the tolerance for the approximation of $\varepsilon$ is $\tau=0.01$.
\end{itemize}

First, let us analyze GDM scenarios with less than 10 DMs involved. In Table \ref{tab:Sim_GDM}, we can see that the error committed when using our methodology is under $0.007$ in comparison to the symmetric (linear) case, whereas it is slightly higher (under $0.01$) in the BLP case. Regarding the average computational time (in milliseconds), our proposal cannot beat linear time, although it is close. However, {\sc{ApOWAMCC}} is much faster than the BLP-based model. Indeed, the elevated computational time required by such a BLP model makes it unfeasible to even test it with more than 10 DMs. For this reason, we have not considered that model in the LSGDM scenario.
\begin{table}[H]
    \centering
    \caption{Results of the simulations for the Linear and BLP versions in small-scale GDM}
    \begin{adjustbox}{max width=\textwidth}
    \begin{tabular}{|c|c|c|c|c|}\hline
        Model &  $n$  & Cost Difference  & {\sc{ApOWAMCC}} Time & Model Time\\\hline
        Linear & 4 & $0.0064 \pm 0.0076$ & $2.18 \pm 4.001$ & $1.26 \pm 3.4484$ \\ 
        Linear & 6 & $0.0058 \pm 0.0053$ & $2.69 \pm 4.165$ & $0.96 \pm 2.785$ \\ 
        BLP & 4 & $0.0085 \pm 0.0104$ & $2.68 \pm 5.2626$ & $3.68 \pm 5.4493$ \\ 
        BLP & 6 & $0.0092 \pm 0.0086$ & $3.76 \pm 5.9052$ & $146.24 \pm 218.6067$ \\ 
    \hline
    \end{tabular}
    \end{adjustbox}
    \label{tab:Sim_GDM}
\end{table}

According to Table \ref{tab:Sim_LSGDM}, the computational complexity of our proposal is nearly the same as the linear model, regardless of the number of DMs, whereas the error in the cost function decreases when the number of variables increases, which makes of {\sc{ApOWAMCC}} a feasible approximation for solving the general version of OWA-MCC in LSGDM scenarios. 

\begin{table}[H]
    \centering
    \caption{Results of the simulations for the Linear version in LSGDM}
    \begin{adjustbox}{max width=\textwidth}
    \begin{tabular}{|c|c|c|c|c|}\hline
        Model &  $n$  & Cost Difference  & {\sc{ApOWAMCC}} Time & Model Time\\\hline
        Linear & 40 & $0.0039 \pm 0.0037$ & $3.02 \pm 3.7605$ & $1.78 \pm 3.2893$ \\ 
Linear & 80 & $0.0021 \pm 0.003$ & $4.43 \pm 6.4858$ & $2.99 \pm 5.3172$ \\ 
Linear & 200 & $0.0013 \pm 0.0024$ & $9.44 \pm 7.1539$ & $10.22 \pm 7.885$ \\ 
Linear & 500 & $0.0002 \pm 0.0002$ & $29.61 \pm 7.125$ & $42.71 \pm 9.824$ \\  
    \hline
    \end{tabular}
    \end{adjustbox}
    \label{tab:Sim_LSGDM}
\end{table}

\section{Conclusions}\label{sec:conclusion}
This study has introduced mutual consensus as a robust consensus measure, positioned as the strongest among classical measures in the literature. Unlike traditional consensus measures that focus on either the distance to the collective opinion or averages of pairwise comparisons, mutual consensus captures the maximum disagreement across opinions, ensuring a more equitable and non-compensatory evaluation of consensus.

Building upon this foundational concept, we have developed new MCC models that integrate mutual consensus as a central constraint. The theoretical analysis revealed key properties of these models, highlighting their potential to enhance decision-making processes.

An interesting application of mutual consensus is the application to approximate the solution of OWA-MCC models. In this regard, we have proposed a linearized version of the OWA-MCC model under symmetry conditions that enhances its computational efficiency. However, we identified cases where linearization is infeasible by showing the inherent complexity and non-convexity of the feasible region in general scenarios.

To address this shortcoming, the relationship between mutual consensus and classical OWA-based MCC models has been analyzed. This theoretical insight guided the development of an MCMC-based algorithm to approximate solutions for the OWA-MCC problem. The algorithm, based on mutual consensus, offers a practical and computationally efficient alternative to classical methods, particularly in LSGDM contexts.

Future research could explore the integration of mutual consensus into other decision-making frameworks, including scenarios with asymmetric cost structures, uncertain preferences, or other complex aggregation operators.

    
\section*{Acknowledgements}
Bapi Dutta acknowledges the support of the Spanish Ministry of Science, Innovation and Universities, and the Spanish State Research Agency through the Ramón y Cajal Research grant (RYC2023-045020-I), Spain.
\bibliographystyle{model2-names}
\bibliography{biblio}
\end{document}